\documentclass[12pt]{article}

\newtheorem{Co}{Conjecture}

\begin{document}

\title{Solution of Hypergraph Turan problem }

\date{July 5,  2015}
\author{Vladimir Blinovsky\footnote{The author was supported by the Sao Paulo Research Foundation (FAPESP), Project  no  2014/23368-6   and
  NUMEC/USP (Project MaCLinC/USP).
}}
\date{\small
 Instituto de Matematica e Estatistica, USP,\\
Rua do Matao 1010, 05508- 090, Sao Paulo, Brazil\\
 Institute for Information Transmission Problems, \\
 B. Karetnyi 19, Moscow, Russia,\\
vblinovs@yandex.ru}

\maketitle\bigskip

\begin{center}
{\bf Abstract}
\bigskip

Using original {\it Symmetrical Smoothing Method} we solve hypergraph $(3,k)$ - Turan problem. \end{center}

Let $X$ be a finite set $|X|=n$. Define ${X\choose k}$ to be the family of $k$ - element subsets of $X$. We say that family (hypergraph) ${\cal A}\subset {X\choose 3}$ satisfies $(3,k)$ - Turan property if  for  the arbitrary
$A\in{[n]\choose k}$ it follows that
$$
|B\in{\cal A}: B\in A |<{k\choose 3}.
$$
Denote the family of hypergraphs which satisfies $(3,k)$ - Turan property by $\hbox{Turan}(3,k)$.
There is a number of sites and conferences devoted to this problem (see~\cite{9},~\cite{10},~\cite{11}).
 We investigate the following
 \bigskip 

{\bf Problem}. For given $n>k>3$ find the maximal (at least one) family which satisfies $(3,k)$ - Turan property.
\bigskip 

This is famous Turan problem.  Actually Turan problem in general case (not only $3$ - hypergraph but $m$ - hypergraph) is the key problem in Erdos's extremal combinatorics. And the case $m=3$ is the first nontrivial
case of great importance. For the surveys and references see~\cite{1},~\cite{2}.

Next we assume that $(k-1)|n$. For other cases it is necessary to find proper  explicit symmetrical constructions
for the optimal Turan hypergraph. I am lazy to do this and it is not needed for asymptotic Turan problem (see~\cite{3},~\cite{4}).

To solve this problem we make some preliminary preparations.

We use the natural bijection between $2^{[n]}$ and $\{ 0,1\}^n$ and don't make difference between these two sets.

Define
$$
\varphi (x)=\frac{1}{\sqrt{2\pi}}\int_{-\infty}^{x/\sigma} e^{-\xi^2 /2}d\xi 
$$
and
\begin{eqnarray}
N(\{\beta_i \} )&=&\sum_{i=1}^{{n\choose 3}} \sum_{x\in {[n]\choose 3}}\varphi\left((x,\beta_i )-3 \right); \label{rt}\\
R(\{ \beta_i \}  ) &=& \sum_{K\in{[n]\choose k}}\varphi\left( \sum_{x\in {K\choose 3}}\varphi\left(\sum_{i=1}^{{n\choose 3}}\varphi\left((x,\beta_i )-3\right)-1/2\right)-1/2\right)-\frac{1}{2}{n\choose k}.\nonumber
\end{eqnarray}
  It can be easily seen
that $N(\{\beta_i\} )$ and $R(\{\beta_i \})$ are convex functions of $\{\beta_{i}\}$. For the arbitrary $\beta_i \in R_+^n,\ \sum_{j=1}^n \beta_{i,j} =3$ hyperspace $\Gamma(\beta_i )=\{y\in R^n\:\ (y,\beta_i )\geq 3\}$  has the property  $\Gamma (\beta_i)\bigcap \{0,1\}^n \leq x$ for some $x\in\{ 0,1\}^n$.  Here the order $\leq$ on subsets $2^{[n]}$ defined by the inclusion.
Condition $R(\{\beta_i\}) =\delta\to 0$ as $ \sigma\to 0$ means that the complement hypergraph $\tilde{T}$ has the property, that there is no $K\in{[n]\choose k}$ s.t. $K\bigcap \tilde{T}=\emptyset$. 
Necessary and sufficient condition for minimization of $N(\{\beta_i\})$ are the Kuhn - Tucker conditions:
\begin{eqnarray}
\label{e1}
&&N^\prime_{\beta_{i,j}}(\{\beta_i\})=\lambda R^\prime_{\beta_{i,j}} ,\ \lambda\in R;\\
&&R(\{\beta_i\})= \delta;\nonumber\\
&&\nonumber \beta_{i,j}\geq 0,\ \sum_{j=1}^n \beta_{i,j}=3,
\end{eqnarray}
where $\delta <0,\ |\delta | <<1$.
If $\{\bar{T}\}=\left\{\arg\max_{T\in\hbox{\hbox{Turan}(3,k)}}{|T|}\right\}$, then for each complement hypergraph $ \tilde{T}={[n]\choose 3}\setminus\bar{T}$  volume $|\tilde{T}|$ is the volume of the  solution of Kuhn - Tucker problem~(\ref{e1}).

Consider the following well known construction of $(3,k)$ - Turan hypergraph. We divide set of vertices $[n]$ into $k-1$ equal parts $B_1 ,\ldots , B_{k-1}$ of size $n/(k-1)$. 
Hypergraph $\bar{T}$ consists of edges in each part $B_r$ and all edges such that each of them has two vertices in $B_r$ and one vertex in $B_{(r+1) |\hbox{mod} (k-1)}$.

It is easy to see, that $\bar{T}$ satisfies $(3,k)$ -Turan property. 

Note, that we find solution of~(\ref{e1})  as the family of $M$ edges $x\in {[n]\choose 3}$, determined by inequalities $(x,\beta_i)\geq  3,\ i\in{n\choose 3}$, where  $\beta_{i,j}\geq 0,\ i\in M$  and $\sum_{j=1}^n \beta_{i,j}=3$.  Number $M$ specifies the number of $\biggl|\left\{i\in {n\choose 3}:\ \exists x\in{[n]\choose 3},\ \hbox{s.t.}\ (x,\beta_i)\geq 3\right\}\biggr|$ .
 
We have
\begin{eqnarray*}&&
N^\prime_{\beta_{i,j}}(\{\beta_i\}) = \frac{1}{\sigma\sqrt{2\pi}}\left(\sum_{x\in{[n]\choose 3},\ x_j =1,\ x_n =0}e^{-\frac{((x,\beta_i)-3)^2}{2\sigma^2}} - \sum_{x\in{[n]\choose 3},\ x_j =0,\ x_n =1}e^{-\frac{((x,\beta_i)-3)^2}{2\sigma^2}}\right);\\
&& R^\prime_{\beta_{i,j}}(\{\beta_i\}) \\&&=\frac{1}{\sigma^3(2\pi)^{3/2}}\Biggl(\sum_{K\in{[n]\choose k}}e^{-\frac{\left(\sum_{x\in{K\choose 3}}\varphi\left(\sum_{i=1}^{{n\choose 3}}\varphi ((x,\beta_i)-3)-1/2\right)-1/2\right)^2}{2\sigma^2}}\\&& \times\Biggl( \sum_{x\in {K\choose 3},\ x_j =1,\ x_n =0}e^{-\frac{\left(\sum_{i=1}^{{n\choose 3}} \varphi((x,\beta_i )-3)-1/2\right)^2}{2\sigma^2}}e^{ -\frac{((x,\beta_i)-3)^2}{2\sigma^2}}\\&&  -\sum_{x\in {K\choose 3},\ x_j =0,\ x_n=1}e^{-\frac{\left(\sum_{i=1}^{{n\choose 3}} \varphi((x,\beta_i )-3)-1/2\right)^2}{2\sigma^2}}e^{ -\frac{((x,\beta_i)-3)^2}{2\sigma^2}}\Biggr)\Biggr).
\end{eqnarray*}
It follows that if $\beta_i =(\beta_{i,1},\ldots ,\beta_{i,n})$ is s.t. $\beta_{i,j}\in \{0,1\}^n,\ \sum_{j=1}^n \beta_{i,j}=3,\ \beta_{i,j_1}=\beta_{i,j_2}=\beta_{i,j_3} =1$ or $\beta_{i,j}=\frac{3}{n},\ j\in [n]$, then there exists $\lambda\in R$, s.t.
\begin{equation}\label{ppp2}
 N^\prime_{\beta_{i,j}}(\{\beta_i\}) =  \lambda R^\prime_{\beta_{i,j}}(\{\beta_i\}),\  i\in {n\choose 3},\ j\in [n-1].  
\end{equation}
Define the set of beta's ${\cal B}= \biggl\{\beta_i =(\beta_{i,1},\ldots ,\beta_{i,n})\in \{0,1\}^n:\ \sum_{j=1}^n \beta_{i,j}=3$ \hbox{and if} $ \beta_{i,j_1}=\beta_{i,j_2}=\beta_{i,j_3} =1,\ j_1,\ j_2 \in B_{r_1},\  j_3 \in B_{r_2},\ j_1 < j_2$ \hbox{following four cases  occur}:\  $ k-2 > r_2 > r_1+1$ or $  r_2 < r_1 -1<k-2$ \hbox{or} $ r_1 = k-1,\ 1<r_2<k-2$
 \hbox{and  third case}  $r_1 <r_2 < r_3$,\ \hbox{ all another beta's,  which don't satisfies these three conditions are zero} $n$-tuples: $\beta_i =\bar{\frac{3}{n}}\biggr\}$. Then
$$
 R(\{\beta_i\})\stackrel{\Delta}{\longrightarrow} 0
 $$
 as $\sigma\to 0$. 
 Hence set ${\cal B}$ is the solution of Kuhn - Tucker problem~(\ref{e1}) and $N(\{\beta_i \})$ achieved global minimum on this set of beta's.

 From here it follows that there exists (unique) $\lambda$ such that equations (\ref{e1}) are satisfied and hence construction of $\tilde{T}$ is optimal.

Easy consequence of optimality of $\bar{T}$ is the validity of the following famous Turan 
\begin{Co}
$$
\lim_{n\to\infty}\frac{|\bar{T}|}{{n\choose 3}}=1-\left(\frac{2}{k-1}\right)^2 .
$$
\end{Co}
At last note that original Symmetrical Smoothing Method allows to prove optimality of the solution of extremal
problems for (binary and not only) sequences in m cases when sufficiently symmetric constructions which are optimal are suggested. We are going to show how to solve many such problems in forthcoming papers.  

    \end{document}